\documentclass[11pt]{article}
\usepackage {graphicx}
\def\ssk{\smallskip}  \def\msk{\medskip}    
\def\text#1{\;\;\hbox{#1}\;\;}    
\def\state #1.{\noindent{\bf#1.\enspace}}     \def\txt#1{\;\hbox{#1}\,}
\def\lset{\big\{\,}    \def\mset{\,\big|\,}   \def\rset{\big\}}

\outer\def\proclaim #1. #2\par{\medbreak
     \noindent{\bf#1.\enspace}{\sl#2}\par
     \ifdim\lastskip<\medskipamount \removelastskip\penalty55\medskip\fi}
\def\paritem#1{\vskip0cm\noindent\hskip20pt{{\rm #1}}\hskip5pt}
\def\eop{\hfill{$\vcenter{\hrule height1pt \hbox{\vrule width1pt height5pt
   \kern5pt \vrule width1pt} \hrule height1pt}$} \medskip}
\def\low#1{{\lower1pt \hbox{$\scriptstyle #1$}}}
\def\Low#1{{\lower2pt \hbox{$\scriptstyle #1$}}}

\def\half{{{}\raise 1pt \hbox{${{\scriptstyle1}\over{\scriptstyle 2}}$}}}
\def\eqalign#1{\begin{array}{lcr} #1 \end{array}}
\def\cC{{\cal C}}     \def\cX{{\cal X}}   \def\cV{{\cal V}}
\def\cW{{\cal W}}   \def\cZ{{\cal Z}}  \def\cU{{\cal U}}   \def\cY{{\cal Y}}
  \def\cP{{\cal P}} 
\def\Bar{\overline}  \def\Hat{\widehat} 
\def\reals{{I\kern-.35em R}}     \def\ball{{I\kern -.35em B}}
\def\mdot{{\kern-.01em\cdot\kern-.02em}}
\def\pls{{\scriptscriptstyle +}} \def\mns{{\scriptscriptstyle -}}
\def\h#1{\hskip #1pt } 
\def\dnto{{\raise 1pt \hbox{$\scriptstyle \,\searrow\,$}}} 
\def\upto{{\raise 1pt \hbox{$\scriptstyle \,\nearrow\,$}}} 
\def\argmin{\mathop{\rm argmin}}

\def\phi{\varphi}  \def\epsilon{\varepsilon}  \def\eps{\varepsilon}         
     \def\qd{\;\;}
\def\tto{\;{\lower 1pt \hbox{$\rightarrow$}}\kern -11pt
           \hbox{\raise 2.8pt \hbox{$\rightarrow$}}\;}
\def\iff{\quad\hbox{$\Longleftrightarrow$}\quad} 
 
\def\implies{\quad\hbox{$\Longrightarrow$}\quad}

    \def\dom{\mathop{\rm dom}\nolimits} 
\def\epi{\mathop{\rm epi}\nolimits}    \def\gph{\mathop{\rm gph}\nolimits} 
\def\iint{\mathop{\rm int}\nolimits}    
    \def\cnv{\mathop{\rm cnv}\nolimits} 
\def\dist{\mathop{\rm dist}\nolimits}   
   
\def\conv{\mathop{\rm conv}\nolimits}   
   
\def\downto{{\raise 1pt \hbox{$\scriptstyle \,\searrow\,$}}} 
\def\for{{\raise1pt \hbox{$\scriptstyle \,|\,$}}} 
\def\ds#1{{\displaystyle #1 }} 
\def\pluss{\hskip1pt \raise1pt\vbox{\hrule width6pt \vskip1pt \hrule width6pt}
 \kern-4pt{\lower1pt\hbox{\vrule height6pt \kern1pt\vrule height6pt}}\hskip5pt}

\def\dfnt{\mathop{\rm def}\nolimits}

\begin{document} \vskip0.5in \centerline{\LARGE\bf 
   VARIATIONAL SUFFICIENCY AND 
}\ssk\centerline{\LARGE\bf  
          SOLUTION STABILITY IN OPTIMIZATION 
} \bigskip \bigskip
\begin{center}
 
\textbf{\textit{Mat\'u\v{s} Benko,}}\footnote{Johann Radon Institute for
 Computational and Applied Mathematics, Linz, Austria; \par\quad
 E-mail: {\it Matus.Benko@oeaw.ac.at}}\quad
\textbf{\textit{R.\ Tyrrell Rockafellar}}\footnote{University of
 Washington, Department of Mathematics, Box 354350, Seattle, WA
 98195-4350; \par \quad  E-mail: {\it rtr@uw.edu}, 
\quad URL: sites.math.washington.edu/$\sim$rtr/mypage.html } 
\bigskip
\end{center} \vskip3cm  

\begin{abstract}
Variational stability, in the sense of local good behavior of optimal 
values and solutions in problems of optimization under shifts in 
parameters, is important not only for validating model robustness
in practical applications but also for confidence of outcomes in the 
design of solution algorithms.  Fundamental results are presented here 
about how such stability relates to a recently developed sufficient 
condition for local optimality called strong variational sufficiency.  

 \end{abstract} \hskip0.5in 

\bigskip

\centerline{ MSC2020 classification 90C31}

\bigskip

\centerline{19 February 2026
}

\vskip60pt
\centerline{\it Dedicated to the memory of H\'edy Attouch}
\newpage

\section{ 
                 Introduction } 

In classical analysis, where problems are so often posed in terms of
equations, the importance of understanding how a solution might depend on
parameters in the equations has long been recognized.  The implicit
function theorem addresses differentiable dependence, while concepts of
well-posedness look to continuous dependence as an essential underpinning
for a good mathematical model in applications.  In variational analysis, as
the extension of classical analysis for handling problems of optimization,
the counterpart is understanding how optimal solutions and optimal values
in a local sense might depend on parameters.  Again, besides being of
interest directly in many applications, this is fundamental to the
development of numerical methodology and whether a problem is sufficiently
well posed in solvability.  Inexactness in computations can often be tied 
to errors in evaluating parameters, which amount to errors in problem
identification, yet it's important for an algorithm still to produce
solutions close to exact ones.  Not only continuous dependence on 
parameters but also Lipschitz continuous dependence comes in then, as
corresponding to solution shifts being in proportion to error sizes.

These motivations have led to decades of research on stability of optimal
solutions and optimal values in which ``stability'' refers to Lipschitz
continuous dependence.  Two papers in particular loom large in the picture
of what will be undertaken here:   Poliquin-Rockafellar \cite{Tilt} for
``tilt'' stability (1998) and Levy-Poliquin-Rockafellar \cite{LocalStability} 
for ``full'' stability (2000).  In both, under assumptions of prox-regularity 
and subdifferential continuity, a condition characterizing the concept was 
identified that utilized coderivatives of associated subgradient mappings.  
Those conditions were tantamount to enhanced second-order sufficient 
conditions for local optimality.  Subsequent contributions on those lines 
in \cite{SecondOrderCalculus}, \cite{FullStability}, \cite{MordSarabi}, 
and elsewhere, 
elaborated the applicability of those coderivative conditions in
particular settings.  Our paper \cite{PrimalDual} paired optimal solutions 
with multiplier vectors in ``primal-dual'' stability.   It also appealed to 
a newer form of second-order optimality condition, called ``strong 
variational sufficiency,'' as a substitute for coderivatives in a central 
role, and that will be carried further here.  Variational sufficiency, 
which originated in \cite{Decoupling} (2019), has been shown to induce a 
kind of local duality \cite{HiddenConvexity} and to support solution 
methods like augmented Lagrangian algorithms in nonconvex optimization 
\cite{ALM}.  In nonlinear programming it's equivalent to the classical 
strong second-order sufficient condition, as established in 
\cite{HiddenConvexity}.

Our main results in this paper will reveal what more needs to be 
combined with strong variational sufficiency to get full stability, and 
on the other hand, the lesser kind of stability that strong variational
sufficiency can guarantee just from a constraint qualification.  
Prox-regularity will be in the forefront, as always, but subdifferential
continuity as a separate assumption will be avoided.  New facts about 
prox-regularity will be brought out which are based on, and add to, the 
connections with variational convexity in the recent papers 
\cite{VariProx}, \cite{ProxRegTests}.\footnote{
   All of these topics in second-order variational analysis were deeply 
   influenced in their evolution by the pioneering work of H.\ Attouch 
   that tied epiconvergence of functions to graphical convergence of their 
   subgradient mappings.}

As in our previous work \cite{PrimalDual}, we adopt a very broad format
for parameterized optimization in finite dimensions where the given problem 
is to
$$ 
    \text{minimize} \phi(x,0) \text{with respect to $x$ }
\leqno\bar\cP$$
for a closed proper function $\phi$ on $\reals^n\times\reals^m$ and the
focusis  on local optimality issues around a particular $\bar x$.  We 
view $\bar P$ as embedded in a parameterized family of problems
$$
     \text{minimize} \phi(x,u) -v\mdot[x-\bar x] \text{with respect to $x$ }
\leqno\cP(v,u)$$ 
as variants, so that 
$$
    \bar\cP =\cP(\bar v,\bar u) \text{for} (\bar v,\bar u)=(0,0). 
$$
We invesigate, for $\delta>0$, the properties of
$$\eqalign{
   M_\delta(v,u) &:=\ds{\argmin_{|x-\bar x|<\delta}}
   \lset\phi(x,u) -v\mdot[x-\bar x]\rset,  \cr
   m_\delta(v,u) &:=\ds{\min_{|x-\bar x|<\delta}}
   \lset\phi(x,u) -v\mdot[x-\bar x]\rset.  }
\eqno(1.1)$$

In the terminology of \cite{Tilt} and \cite{LocalStability},  {\it
stability} of $\bar x$ refers to $M_\delta(0,u)$ and $m_\delta(0,u)$ being,
for small enough $\delta$, Lipschitz continuous in $u$ in a neighborhood
of $u=0$, with $M_\delta$ single-valued and given at $0$ by $\bar x$, while 
{\it tilt stability\/} refers to those properties holding for $M_\delta(v,0)$ 
and $m_\delta(v,0)$ for $v$ in a neighborhood of $v=0$;  {\it full 
stability\/} joins these by requiring Lipschitz continuity of 
$M_\delta(v,u)$ and $m_\delta(v,u)$ in $(v,u)$ around 
$(v,u)=(0,0)$.  It will be useful to speak of {\it substability\/} 
when $M_\delta(0,u)$ and $m_\delta(0,u)$ are just continuous in $u$, not
necessarily Lipschitz continuous.  A hybrid soon to come up is
$$\eqalign{
  &\hbox{{\it full substability\/} of $\bar x$:  having 
           $M_\delta(v,u)$ and $m_\delta(v,u)$ just continuous} \cr
  &\h{5}\hbox{in $(v,u)$ but Lipschitz continuous in $v$,
          locally uniformy in $u$. } }
\eqno(1.2)$$
In all cases, $\bar x$ is obviously an {\it isolated\/} locally optimal 
solution, and this persists parametrically.

The challenge is to tie these properties to some condition or other
on the subgradients of $\phi$ in the general sense of variational
analysis \cite{VA}.  In this, the mappings $M$ and $Y$ defined by
$$
   M(v,u) =\lset x \mset Y(x,u,v)\neq\emptyset \rset, \qquad
   Y(x,u,v):= \lset y \mset (v,y)\in\partial\phi(x,u) \rset 
\eqno(1.3)$$
have a big role around
$$
       \Bar Y = Y(\bar x,0,0),
    \text{the set of {\it multiplier vectors\/} $\bar y$ associated with 
                $\bar x$ in $\bar P$.}
\eqno(1.4)$$

The necessary condition $0\in \partial_x\phi(\bar x,0)$ for the local 
optimality of $\bar x$ in $\bar P$ yields a $\bar y\in\Bar Y$ through 
\cite[10.6]{VA} under the {\it basic constraint qualification\/} on 
horizon subgradients, namely
$$
      (0,y) \in \partial^\infty\phi(\bar x,0) \implies y=0.
\eqno(1.5)$$
Then, in fact, the multiplier set $\Bar Y$ is not just nonempty, but also
compact, and useful properties of the set-valued mapping $Y$ for $(x,u,v)$
in a neighborhood of $(\bar x,0,0)$ are made available as well.

The {\it variational sufficient condition\/} for the local optimality of
$\bar x$ holds for a multiplier $\bar y\in\Bar Y$ if there exists 
$e\geq 0$ (``elicitation parameter'') such that the function
$\phi_e(x,u)=\phi(x,u)+\frac{e}{2}|u|^2$ is variationally convex at
$(\bar x,0)$ for the subgradient $(0,\bar y)\in\partial\phi_e(\bar x,0)
=\partial\phi(\bar x,0)$.  The {\it strong\/} version has variational
{\it strong\/} convexity of $\phi_e$.  These sufficient conditions will be 
reviewed in Section 2 along with the details of prox-regularity and 
variational convexity.   They were introduced in \cite{Decoupling} and 
broadly explored in \cite{HiddenConvexity} in the framework of 
{\it generalized nonlinear programming}, where
$$
   \phi(x,u)=f_0(x)+g(F(x)+u) \text{with $g$ closed proper convex,
      $f_0$ and $F$ in $\cC^2$.}
\eqno(1.6)$$ 
Classical nonlinear programming is the case of this in which $g$ is the
indicator $\delta_K$ of the cone $K=\reals^s_\mns\times\reals^{m-s}$.
There, the basic constraint qualification is equivalent to the
Mangasarian-Fromovitz constraint qualification, and the strong variational
sufficient condition is equivalent to the standard strong second-order 
sufficient condition for local optimality \cite[Example 1]{HiddenConvexity}.  
But many other examples of variational sufficiency were developed in
\cite{HiddenConvexity} as well.

We can now state our main results about stability of a locally optimal
solution $\bar x$ to problem $\bar P$ in relation to variational 
sufficiency.  It will be convenient in this to adopt the terminology that
$$\eqalign{
  \text{``variational sufficiency at $\bar x$'' means having 
      it at $\bar x$ for all $\bar y\in\Bar Y$, i.e., 
      $(0,\bar y)\in\partial\phi(\bar x,0)$, }  &\cr
  \text{``full multiplier prox-regularity at $\bar x$'' means 
    $\phi$ prox-regular at $(\bar x,0)$ for all such $(0,\bar y)$.} &}
\eqno(1.7)$$

\proclaim Theorem 1.1 {\rm (variational sufficiency versus full stability
and full substability)}.
Under the basic constraint qualification (1.5),
$$\eqalign{
   \text{(A)\qd full stability holds at $\bar x$ with full
           multiplier prox-regularity} \implies &\cr 
   \text{(B)\qd strong variational sufficiency holds at $\bar x$ } 
             \implies &\cr
   \text{(C)\qd full substability holds at $\bar x$ with full
           multiplier prox-regularity.}  &}
$$
Under (C), hence also (A) and (B), an open convex neighborhood 
$\cV\times\cU$ of $(0,0)$ exists such that, for small enough $\delta$ 
and low enough $\alpha > m_\delta(0,0)=\phi(\bar x,0)$, 
$$
    M_\delta (v,u) = 
   M(v,u)\cap\lset x \mset |x-\bar x|<\delta, \phi(x,u)<\alpha\rset
      =: M^{\delta,\alpha}(v,u),
\eqno(1.8)$$
and $m_\delta(v,u)$ is $\cC^1$ in $v\in\cV$ for each $u\in\cU$
with $\nabla_v m_\delta(v,u)=M_\delta(v,u)$.  Under (B), furthermore 
$m_\delta(v,u)$ is uniformly hypo-convex in $u\in\cU$ for $v\in\cV$ with
$$
  \partial_u m_\delta(v,u)= Y(x,u,v) \text{for} x= M_\delta(v,u),
\eqno(1.9)$$
where the uniform hypo-convexity refers to the existence of $e\geq 0$ 
such that $\,m_\delta(v,u)+\frac{e}{2}|u|^2\,$ is convex in $u\in\cU$ for 
each $v\in\cV$.

The hypo-convexity observation at the end of Theorem 1.1 is new even for 
the case of $\phi$ in (1.6) that corresponds to classical nonlinear 
programming.

A counterexample will show that (B) $\not\Rightarrow$ (A), but whether 
(C) $\Rightarrow$ (B), making that pair equivalent, is open to
serious conjecture.  Another unsettled question is whether the basic 
constraint qualification might automatically hold under full substability,
hence also full stabilty.   It holds under the latter at least for $\phi$
in the composite format (1.6); see Mordukhovich-Sarabi 
\cite[Theorem 5.1]{MordSarabi2} (note that the relevant part of 
\cite[Theorem 5.1]{MordSarabi2} doesn't rely on the special case of $g$ 
considered in the paper).  The basic constraint qualification (1.5) is 
known to be equivalent  the Aubin property of the mapping 
$u \mapsto \epi\phi(\cdot,u)$ at $\bar u$ for the image point 
$(\bar x,\phi(\bar x,\bar u))$ \cite[10.16]{VA}, and in that way connects 
naturally to the behavior of $M_\delta(v,u)$ with respect to $u$.

Since (B) and (C) definitely fall short of (A), there is the question of 
what needs to be combined with them to bring about an equivalence, other 
than an appeal back to the original coderivative conditions in 
\cite{LocalStability}.  We'll provide an answer which, for its 
formulation, requires graphical derivatives.  Recall that the graphical 
derivative of the single-valued truncated mapping $M^{\delta,\alpha}$ 
in (1.8) at $(v,u)$ is the generally set-valued mapping 
$DM^{\delta,\alpha}(v,u): \reals^n\times\reals^m\tto\reals^n$ having, as 
its graph, the tangent cone to the graph of $M^{\delta,\alpha}$ at 
$(v,u,x)$ for $x=M^{\delta,\alpha}(v,u)$.  We will be looking at, more 
specially, the graphical derivatives $DM^{\delta,\alpha}_v(u)$ of the 
``partial mappings''
$$ 
       M^{\delta,\alpha}_v:=M^{\delta,\alpha}(v,\cdot):   
                       u\mapsto x=M^{\delta,\alpha}(v,u),
\eqno(1.10)$$
which have inner and outer norms, given in the notation
$\xi\in DM^{\delta,\alpha}_v(u)(\omega)$ by
$$
  |DM^{\delta,\alpha}_v(u)|^\mns =\ds{\sup_{|\omega|\leq 1}\,
         \inf_\low{\low{\xi\in M^{\delta,\alpha}_v(u)(\omega)}}}|\xi|,  
      \qquad
  |DM^{\delta,\alpha}_v(u)|^\pls =\ds{\sup_{|\omega|\leq 1}\,
              \sup_{\xi\in M^{\delta,\alpha}_v(u)(\omega)}}|\xi|.
\eqno(1.11)$$

\proclaim Theorem 1.2 {\rm (full stability from strong variational 
sufficiency or just full substability)}.
Condition (A) in Theorem 1.1 is equivalent to the combination of 
condition (C), or (B), with
$$
    \limsup_{(v,u)\to (0,0)} |DM^{\delta,\alpha}_v(u)|^\mns <\infty.
\eqno(1.12)$$ 

The proofs of these theorems will be presented in Section 3 after 
groundwork in Section 2 on understanding prox-regularity and its 
connections to variational convexity and tilt stability.

\section{
        Prox-regularity background with new developments
}

Prox-regularity is a primal-dual-localized property of a function that
has received many years of close attention since being introduced in
\cite{ProxReg}, especially for its role in second-order variational
analysis and its connection to behaviors tied, one way or another, to
convexity.  It's connection to {\it variational\/} convexity was 
strengthened in \cite{VariProx} and \cite{ProxRegTests}.  Consequences 
will be taken farther in this section in building a platform for the 
proofs of our main results.

For now, we consider just an extended-real-valued function $f$ on 
$\reals^n$ instead of $\phi$ on $\reals^n\times\reals^m$, supposing it to 
be closed (lower semicontinuous) and proper.  We draw on the theory in
\cite{VA} of general (limiting) subgradients $v\in\partial f(x)$ and 
horizon subgradients $v\in\partial^\infty f(x)$.   Prox-regularity of $f$ at 
$\bar x$ with respect to a subgradient $\bar v\in\partial f(\bar x)$
involves open neighborhoods $\cX$ of $\bar x$ and $\cV$ of $\bar v$
together with an $\alpha>f(\bar x)$ that truncates $\cX$ to
$$
    \cX^\alpha_f = \cX\cap\lset x \mset f(x)<\alpha\rset,   
\eqno(2.1)$$
a so-called  {\it $f$-attentive neighborhood\/} of $\bar x$.  Such 
truncation is important in the absence of continuity of $f$ at $\bar x$ in 
ensuring that localization only captures geometric properties of $\epi f$ 
in a neighborhood of $(\bar x,f(\bar x))$.   In the case of prox-regularity, 
the property is the existence of $r\in(-\infty,\infty)$ such that
$$
  f(x')\geq f(x)+v\mdot[x'-x]-\frac{r}{2}|x'-x|^2, \;
 \forall x'\in\cX, \; (x,v)\in(\cX^\alpha_f\times\cV)\cap \gph\partial f.
\eqno(2.2)$$
More specifically, this is prox-regularity at level $r$.  But $r$ could be
replaced in (2.1) by $s=-r$. in which case this would be the property of
uniform quadratic growth at level $s$, which like $r$ can be negative as
well as positive.  

For $r=s=0$, there is a close resemblance to the global subgradient 
inequality that signals convexity, but in a localization that is
primal-dual and $f$-attentive.  More broadly, the global $s$ version 
with $s>0$ would signal strong convexity at that level, but even for
negative $s$, we can think of it in connection with $f$ being $s$-convex
in the sense of $f(x)-\frac{s}{2}|x|^2$ being convex.  These insights have
motivated the following definition:  $f$ is {\it variationally $s$-convex\/} 
at $\bar x$ for $\bar v\in \partial f(\bar x)$ if, for an $f$-attentive
neighborhood as in prox-regularity, with $\cX$ and $\cV$ open convex,
there is a truly $s$-convex closed proper function $\Hat f$ such that
$$\eqalign{
  \Hat f\leq f \text{on $\cX$,} \qd (\cX\times\cV)\cap \gph\partial\Hat f = 
          (\cX^\alpha_f\times\cV)\cap \gph\partial f &\cr
 \hbox{and, for all $(x,v)$ in this truncated graph, $\Hat f(x)= f(x)$.}  &}
\eqno(2.3)$$ 

This concept surfaced already in 1998 in \cite{Tilt} for $s>0$, but didn't
get a name until 2019 in \cite{VarConv}.  Full articulation came in
\cite{VariProx}, where it was coordinated with $\partial f$ being {\it
maximally $s$-monotone\/} locally at $\bar x$ for 
$\bar v\in\partial f(\bar x)$ in the sense that, for some $f$-attentive
neighborhood $\cX_f^\alpha\times\cV$ of $(\bar x,\bar v)$ with open $\cX$
and $\cV$,
$$
 (v'-v)\mdot(x'-x)\geq s|x'-x|^2 \text{when} (x,v),\,(x',v')
               \in(\cX_f^\alpha\times\cV)\cap\gph\partial f,
\eqno(2.4)$$
but without this holding for any strictly larger subset of $\cX\times\cV$.

\proclaim Theorem 2.1 {\rm (equal reflections of prox-regularity 
\cite{VariProx})}.
For any $r\in\reals$ and $s=-r$, the following properties at $\bar x$ for
$\bar v\in\partial f(\bar x)$ are equivalent:
\paritem{(a)} prox-regularity of $f$ at level $r$,
\paritem{(b)} variational convexity of $f$ at level $s$,
\paritem{(c)} $f$-attentive locally maximal $s$-monotonicity 
               of $\partial f$.

It's {\it not\/} claimed that (a), (b) and (c) are always equivalent for 
the same $\cX^\alpha_f\times\cV$.  Shifts may be needed in passing between 
the three properties, but there is nonetheless out of the equivalence a 
single value $\cnv f(\bar x\for\bar v) \in [\infty,\infty]$ such that
$$
   \cnv f(\bar x\for\bar v)=-\left\{\h{-17}\eqalign{
        &\hbox{liminf of the available $r$-values in (2.2) as} \cr
        &\hbox{$\;\cX\times\cV$ shrinks down to $(\bar x,\bar v)$ 
          and $\alpha\dnto f(\bar x)$,} }\right.
\eqno(2.5)$$ 
$$
   \cnv f(\bar x\for\bar v)=\;\left\{\h{-17}\eqalign{
        &\hbox{limsup of the available $s$-values in (2.3) as} \cr
        &\hbox{$\;\cX\times\cV$ shrinks down to $(\bar x,\bar v)$ 
          and $\alpha\dnto f(\bar x)$,} }\right.
\eqno(2.6)$$
$$
   \cnv f(\bar x\for\bar v)=\;\left\{\h{-17}\eqalign{
        &\hbox{limsup of the available $s$-values in (2.4) as} \cr
        &\hbox{$\;\cX\times\cV$ shrinks down to $(\bar x,\bar v)$ 
          and $\alpha\dnto f(\bar x)$.} }\right.
\eqno(2.7)$$
This common value $\,\cnv f(\bar x\for\bar v)$ is the {\it convexity 
modulus\/} of $f$ at $\bar x$ for the subgradient  $\bar v$.  There are 
more equivalences in the theory in \cite{VariProx} beyond those in 
Theorem 2.1, but this is all we'll need here.

A distinction in terminology requires care, however.  Saying that $f$ 
is prox-regular at $\bar x$ for a subgradient $\bar v\in\partial f(\bar x)$ 
means prox-regularity holds there at  {\it some\/} level $r\in\reals$, but 
saying that $f$ is variationally convex at $\bar x$ for $\bar v$
specifically means variational convexity at level $s=0$.   Level $s>0$ is 
variational {\it strong\/} convexity.  Level $s<0$ is variational 
{\it hypo\/}-convexity.  By itself without a designated function and 
location, ``variational convexity'' refers instead to the general topic at 
all levels, which through Theorem 2.1 is identical to the general topic of
prox-regularity.

Variational strong convexity of $f$ at $\bar x$ for $\bar v=0$ as
subgradient is important for understanding tilt stability of $f$ at 
$\bar x$ in the sense of the mappings 
$$\eqalign{
   m_\delta(v) &:=\ds{\inf_{|x-\bar x|<\delta}}
   \lset f(x) -v\mdot[x-\bar x]\,\rset,  \cr
   M_\delta(v) &:=\ds{\argmin_{|x-\bar x|<\delta}}
   \lset f(x) -v\mdot[x-\bar x]\,\rset,  }
\eqno(2.8)$$  
being, for small enough $\delta>0$, Lipschitz continuous on a neighborhood
of 0 with $M_\delta$ single-valued and assigning $\bar x$ to 0.
Originally in \cite[Theorem 1.3]{Tilt}, variational strong convexity,
without it yet having that name, was shown to be necessary and sufficient
for tilt stability under the assumption not only of prox-regularity of $f$
at $\bar x$ for $\bar v=0$ but also {\it subdifferential continuity\/}
there, which refers to having
$$
    f(x)\to f(\bar x) \text{when} (x,v)\to (\bar x,\bar v)
             \text{within} \gph\partial f.
\eqno(2.9)$$
The additional assumption has turned out to be superfluous, however.

\proclaim Theorem 2.2 {\rm (tilt stability without assuming subdifferential 
continuity, Gfrerer \cite[Theorem 5.1]{Gfrerer})}.
The variational strong convexity of $f$ at $\bar x$ for subgradient
$\bar v=0$ is equivalent to $f$ being prox-regular at $\bar x$ for 
$\bar v=0$ and tilt stable there.  In more detail, this means having an 
interval $(0,\sigma)$ such that variational $s$-convexity of $f$ for all 
$s\in (0,\sigma)$ corresponds to $M_\delta$ being, for all 
$s\in (0,\sigma)$, Lipschitz continuous with constant $s^{-1}$ in 
some neighborhood of 0.

This result relies on facts about prox-regularity that concern a modified
version of the subdifferential continuity in (2.9), namely
$$\eqalign{
  \text{under the prox-regularity condition (2.2),} &\cr
\quad\text{(a)\qd $f(x)$ depends continously on} 
              (x,v)\in (\cX_f^\alpha\times\cV)\cap\gph\partial f,  &\cr
\quad\text{(b)\qd $(\cX_f^\alpha\times\cV)\cap\gph\partial f$ 
        is closed relative to $\cX\times\cV$,  } &}
\eqno(2.10)$$
where (a) comes from \cite[Proposition 2.3]{ProxReg} and (b) has been
added by Gfrerer in \cite[Lemma 3.4]{Gfrerer}.  (In the cited works, $\cX$ 
and $\cV$ are open $\eps$-balls around $\bar x$ and $\bar v$, and 
$\alpha=f(\bar x)+\eps$, but the broader formulation in (2.10) is supported
by the simple observation in Proposition 2.5(a) below.)  In particular, 
(2.10) covers the version of (2.9) restricting to $f(x)<\alpha$.  Another 
consequence is part (b) of the theorem coming next, which records additional 
facts needed later.  

\proclaim Theorem 2.3 {\rm (properties joined to tilt stability)}.
  \paritem{(a)} The function $m_\delta$ in (2.8) is always finite, concave
and globally Lipschitz continuous on $\reals^n$.  In tilt stability, it is
moreover differentiable on a neighborhood
of 0 with
$\nabla m_\delta(v) = -M_\delta(v)$.
  \paritem{(b)} Under prox-regularity of $f$ at $\bar x$ for subgradient
$\bar v=0$, there exists $\alpha > m_\delta(0)=f(\bar x)$ such that,
in tilt stability with $v$ close enough to $0$, 
$$
    M_\delta(v)= (\partial f)^{-1}(v)\cap
                 \lset x \mset |x-\bar x|<\delta, \; f(x)<\alpha\rset. 
\eqno(2.11)$$
When $f$ is subdifferentially continuous at $\bar x$, this holds with the 
$\alpha$ replaceable by $\infty$. 
  \paritem{(c)} Variational strong convexity at level $s>0$ corresponds in
the equivalence in Theorem 2.2 to having $s^{-1}$ be a Lipschitz constant for
$M_\delta$ when $\delta$ is sufficiently small.

\state Proof.  (a):  The concavity is evident because the formula presents
$m_\delta$ as the pointwise infimum of a collection of affine functions.
A lower semicontinuous proper function is bounded from below on every 
compact subset, hence $f$ is bounded from below by some $\beta$ on the ball in
(2.8).  In that case,
$$
    \infty>f(\bar x)\geq m_\delta(v)\geq \beta -\delta|v|,
$$
which tells us that $-m_\delta$ is a finite convex function on $\reals^n$
with recession function $\leq |\cdot|$.  Hence
$$ 
  -m_\delta(v') \leq -m_\delta(v) +\delta|v'-v|
                 \text{for all} v,\,v',
$$
signaling Lipschitz continuity with constant $\delta$.  That property
of $-m_\delta$ carries over to $m_\delta$.   

For the claim in (a) about $M_\delta(v)$, we note that if $x$ is a
minimizer for $v$ in (2.8), $x\in M_\delta(v)$, then 
$m_\delta(v) = f(x)-v\mdot[x-\bar x]$, while for other $v'$ merely 
$m_\delta(v')\leq f(x)-v'\mdot[x-\bar x]$.  That combination implies 
$m_\delta(v')-m_\delta(v)\leq [v-v']\mdot x$, according to which $x$ is a 
subgradient of the convex function $-m_\delta$ at $v$.  The set of 
subgradients is convex, so that if $M_\delta(v)$ consists of a unique $x$,
as in tilt stability, that $x$ is the unique subgradient of
$-m_\delta$ at $v$.  In convex analysis this corresponds to $-m_\delta$
being differentiable at $v$ with $x$ as its gradient there.  Then the
vector $-x = -M_\delta(v)$ is the gradient of $m_\delta$ at $v$.

(b):   This takes advantage of the equivalence in Theorem 2.2, which
provides a strongly convex function $\Hat f$ satisfying (2.3) for $\bar v=0$
and some $\cX$, $\cV$ and $\alpha > f(\bar x)$ as described.  Then 
$f(\bar x)\geq m_\delta(0)$, with equality holding under tilt stability.   
From part (a), we get $m_\delta(v)<\alpha$ for $v$ in a neighborhood of 
$0$ inside $\cV$, so that having $x\in M_\delta(v)$ for $\delta$ such
that $\lset x\mset |x-\bar x|<\delta\rset\subset\cX$ will entail having 
$v\in\partial f(x)$ with $(x,v)\in\cX_f^\alpha\times\cV$.  That's the same 
in (2.3) as having $v\in\partial\Hat f(x)$ with $(x,v)\in\cX\times\cV$.  
The convexity of $\partial\Hat f$ associates that in turn with $x$ giving 
the minimum of $\Hat f(x)-v\mdot[x-\bar x]$ over $\cX$ and hence, by the 
relationship in (2.3), also the minimum of $f(x)-v\mdot[x-\bar x]$ over $\cX$.  That leads 
to the equation in (2.11).  In the presence of the subdifferential 
continuity in (2.9), just having $(x,v)\in\gph\partial f$ near enough to 
$(\bar x,0)$ ensures having $f(x)<\alpha$.  The $\alpha$ condition in (2.11)
 is then satisfied and doesn't need to be entered explicitly.

(c):  Again this utilizes the equivalence in Theorem 2.2, with the
function $\Hat f$ in (2.3) being $s$-convex.  The $s$-convexity with $s>0$
implies
$$
 (v'-v)\mdot(x'-x)\geq s|x'-x|^2 \text{when} (x,v),\,(x',v')
               \in(\cX\times\cV)\cap\gph\partial\Hat f,
$$
with the consequence that $|x'-x|\leq s^{-1}|v'-v|$.  Through (2.3) and
part (b), this carries over in (2.3) to $x=M_\delta(v)$ and 
$x'=M_\delta(v')$.    \eop 

In connection with Theorem 2.3(c), we note that this property doesn't
quite follow from Gfrerer's equivalence in Theorem 2.2, because that
applies at a level $s$ only when variational strong convexity of $f$ also 
holds at some level $s'>s$.

The original characterization of tilt stability in \cite[Theorem 1.3]{Tilt},
which utilized subdifferential continuity, included also equivalence with 
the positive-definiteness of the Mordukhovich second-order subdifferential 
$\partial^2 f(\bar x\for 0)$.  A substitute for that subdifferential is 
needed in the absence of subdifferential continuity, because the graph of 
$\partial^2 f(\bar x\for 0)$, being defined in terms of limits of normal 
cones to the graph of $\partial f$ that aren't $f$-attentive, can be much 
larger than the graph of $\partial^2 f^\alpha(\bar x\for 0)$ for 
$f^\alpha(x)=f(x)$ when $f(x)\leq\alpha$, otherwise $f^\alpha(x)=\infty$.
The remedy is replacing $\partial^2 f(\bar x\for 0)$ by the Gfrerer 
second-order subdifferential, defined in \cite{Gfrerer}, which we'll denote 
here by $\partial_G^2 f(\bar x\for 0)$.  It simply modifies the Mordukhovich 
definition by requiring $f$-attentivity, so that
$
    \partial_G^2 f(\bar x\for 0) = \partial_G^2 f^\alpha(\bar x\for 0) = 
                       \partial^2 f^\alpha(\bar x\for 0).
$
That way, in terms of $\partial_G^2 f(\bar x\for 0)$, the second-order 
subdifferential characterization in \cite{Tilt} carries over without the 
assumption of subdifferential continuity; cf.\ \cite[Theorem 5.1]{Gfrerer}.

Here, however, instead such of coderivative-based subdifferentials we'll
employ for our purposes the {\it strict\/} second-order subdifferential 
$\partial^2_*f(\bar x\for 0)$ that was introduced in \cite{ProxRegTests} by
$$\eqalign{
    (\xi,\mu) \in\gph \partial^2_*f(\bar x\for\bar v) \iff \exists\,
    (x^\nu,v^\nu)\to(\bar x,\bar v) \text{in}\gph\partial f \text{and} &\cr
    (\xi^\nu,\mu^\nu)\to(\xi,\mu),\; \tau^\nu\dnto 0, \text{with}
    (x^\nu+\tau^\nu\xi^\nu,v^\nu+\tau^\nu\mu^\nu)\in\gph\partial f, &\cr
    \hbox{and furthermore both $\,f(x^\nu)\to f(\bar x)$ and
     $\,f(x^\nu+\tau^\nu\xi^\nu)\to f(\bar x)$.} &}
\eqno(2.12)$$
This graphical prescription means that the difference quotient mapping 
$\partial^2_*f(\bar x\for\bar v)$ is the outer graphical limit of the
first-order difference quotient mappings $\Delta_\tau[\partial f](x\for v)$
as $\tau\downto 0$ while also $(x,v)\to (\bar x,\bar v)$ $f$-attentively.
The {\it definiteness\/} modulus of $\partial^2_*f(\bar x\for\bar v)$ is
$$
   \dfnt[\partial^2_*f(\bar x\for\bar v)] = \sup\lset s\in\reals \mset
   \mu\mdot\xi\geq s|\xi|^2 \txt{when}
  \mu\in\partial^2_*f(\bar x\for \bar v)(\xi) \rset.
\eqno(2.13)$$

\proclaim Theorem 2.4 {\rm (alternative criterion for tilt stability)}.
The following are equivalent when $f$ is prox-regular at $\bar x$ with
respect to having 0 as a subgradient there:
  \paritem{(a)} $f$ is tilt-stable at $\bar x$.
  \paritem{(b)} $f$ has a local minimum at $\bar x$, and
               $0\in\partial^2_*f(\bar x\for 0)(\xi)$ only for $\xi=0$.
  \paritem{(c)} $\partial^2_*f(\bar x\for 0)$ is positive-definite:
              $\dfnt[\partial^2_*f(\bar x\for 0)]>0$. 
\newline
In these circumstances, for any positive 
$s<\dfnt[\partial^2_*f(\bar x\for 0)]$ in (c), $s^{-1}$ serves as a 
local Lipschitz constant for the tilt stability in (a). 

\state Proof.  The equivalence of (c) with the variational strong
convexity of $f$ at $\bar x$ for subgradient $\bar v=0$ is known from
\cite[Theorem 3.6]{ProxRegTests} and the fact that prox-regularity makes 
subgradients be proximal subgradients.  In that equivalence, the
variational strong convexity holds at level $s$ for any positive 
$s<\dfnt[\partial^2_*f(\bar x\for\bar v)]$.  Then from Theorem 2.2 we have
(c) being equivalent to (a) with $s^{-1}$ serving as a Lipschitz constant 
locally for $M_\delta$ in the tilt stability.  We also have (b) implied by
the minimization in (a) and the positive-definiteness in (c).  All that 
remains is demonstrating that (b) implies (a).

The assumption that $\bar x$ is a local minimizer of $f$ guarantees that
$\bar x\in M_\delta(0)$ when $\delta$ is small enough.  The nonsingularity 
condition in (b) is equivalent by \cite[Proposition 3.3]{ProxRegTests} to 
the existence of a neighborhood $\cV$ of 0 and an $f$-attentive 
neighborhood $\cX$ of $\bar x$ such that the mapping 
$M^*: v\in \cV \mapsto \cX_f^\alpha\cap(\partial f)^{-1}(v)$ 
is single-valued and Lipschitz continuous.  Since 
$m_\delta(v) \leq f(\bar x)<\alpha$ for all $v$, we will therefore have, 
when $\delta$ is small enough, that $x\in M_\delta(v)$ and $v\in\cV$ 
entail $x\in\cX_f^\alpha$ along with $v\in \partial f(x)$, hence 
$x=M^*(v)$.  This confirms that {\it to the extent that $M_\delta$ is
nonempty-valued\/} it is single-valued and Lipschitz continuous, 
assigning $\bar x$ to 0.  This leaves us with demonstrating that the
nonemptiness can be arranged by taking $\delta$ still smaller, if necessary.

From $\bar x$ being the unique minimizer of $f$ over the open
$\delta$-ball around $\bar x$, we know that $f(x)>f(\bar x)$ when 
$|x-\bar x|<\delta$.  Taking $\delta$ to be slightly smaller, we can
arrange that $\bar x$ is the unique element of $\argmin f_\delta$, with
$f_\delta$ being the sum of $f$ and the indicator of the closed
$\delta$-ball around $\bar x$.  For any $v$, the function 
$f_\delta^v(x)=f_\delta(x)-v\mdot[x-\bar x]$ is lower semicontinuous with
bounded domain, so that $\argmin f_\delta^v \neq\emptyset$.   Moreover, 
$\argmin f_\delta^v = M_\delta(v)$ when 
$\argmin f_\delta^v \subset \lset x \mset x-\bar x|<\delta \rset$.  Our goal
can be achieved by verifying that the latter must hold for sufficiently
small $\delta$.   The key to that is the fact that $f_\delta^v$ converges
epigraphically to $f_\delta$ as $v\to 0$ by \cite[7.8a]{VA}.  This ensures
by \cite[7.33]{VA} that the outer limit of the nonempty compact sets 
$\argmin f_\delta^v$ as $v\to 0$ lies within $\argmin f_\delta$, which is
just $\{\bar x\}$.   \eop

The remainder of this section turns from tilt stability to more basic
consequences of prox-regularity and their impact on subgradient calculus.
We note first that, although the definitions of levels of prox-regularity 
and variational convexity center on a particular pair $(\bar x,\bar v)$ 
in $\gph\partial f$, with $\cX_f^\alpha\times\cV$ being a sufficiently 
small neighborhood, the properties in (2.2) and (2.3) have a life beyond 
just that, because they can operate for multiple pairs simultaneously.

\proclaim Proposition 2.5 {\rm (simultaneous prox-regularity)}.
Let the prox-regularity condition (2.2) hold for some open set 
$\cX\times\cV$ and $\alpha\in\reals$.  Then for every 
$\bar x\in\cX_f^\alpha$,
 \paritem{(a)} $f$ is prox-regular at $\bar x$ for every 
$\bar v\in\cV\cap\partial f(\bar x)$,
 \paritem{(b)} $\partial f(\bar x)$ is convex if 
$\cV\supset\partial f(\bar x)$, while $\cV\cap\partial f(\bar x)$ is 
convex if $\cV$ is convex.

\state Proof.   (a):  This is obvious from the fact that $\cX_f^\alpha$ 
is an $f$-attentive neighborhood of each of its elements.  (b):   Any 
vector in the convex hull of $\cV\cap\partial f(\bar x)$ is a convex 
combination $\bar v=\sum_{i=0}^n\lambda_i \bar v_i$ of vectors 
$\bar v_i\in\cV\cap\partial f(\bar x)$.  By taking the same convex 
combination of the inequalities (2.2) holding for each $i$, we see the 
inequality holding for $\bar v$.  That resultant inequality says in 
particular that $\bar v$ is a proximal subgradient, so in particular
$\bar v\in\partial f(\bar x)$.  When $\cV$ is convex, $\bar v$ belongs 
also to $\cV$, but even when $\cV$ isn't convex, if 
$\cV\supset\partial f(\bar x)$ we anyway have 
$\cV\cap\partial f(\bar x)= \partial f(\bar x)$ and can conclude that
the convex hull of $\partial f(\bar x)$ is a subset of $\partial f(\bar x)$.
That means $\partial f(\bar x)$ is convex.   \eop

Analogous statements hold also for the variational convexity condition 
(2.3) and for the local maximality in the monotonicity condition (2.4).  
The two counterparts to (a) are just as obvious.  The counterparts for 
(b) are clear from the convexity of $\partial\Hat f(x)$ in (2.3) and the 
maximality imposed in (2.4), since that inequality would persist in 
adding any $v\in\cV$ that was missing from $\cV\cap\partial f(\bar x)$ 
but belonged to its convex hull.  

The statement and proof of the next result make use of the convexity 
modulus $\cnv f(\bar x\for\bar v)$ given simultaneously by (2.5), (2.6)
and (2.7).

\proclaim Proposition 2.6 {\rm (extended prox-regularity)}.
If $f$ is prox-regular at $\bar x$ for all $\bar v\in C$  for some compact 
set $C\subset\partial f(\bar x)$, or equivalently 
$\cnv f(\bar x\for\bar v)>-\infty$ for such $\bar v$, then 
 \paritem{(a)} $\conv C\subset\partial f(\bar x)$ and $f$ is prox-regular 
at $\bar x$ for all $\bar v \in \conv C$.   
 \paritem{(b)} $\exists\,\cX$, $\cV$, $\alpha$ and $r$ for which 
the prox-regularity condition (2.2) holds with $\cV\supset\conv C$.
 \paritem{(c)} Here $\cX$ and $\cV$ can be taken to be convex and,
together with $\alpha$ and $r$, adjusted so that the variational convexity 
condition (2.3) for $s=-r$ holds along with (2.2) for all 
$\bar v\in\conv C$, moreover with $s >\bar s$ if 
$\cnv f(\bar x\for\bar v)>\bar s$ for all $\bar v\in C$.

\state Proof.   (a):  Extending in the argument for Proposition 2.5(b), 
we represent an arbitrary $\bar v\in\conv C$ as a convex combination 
$\sum_{i=0}^n\lambda_i\bar v_i$ of vectors $\bar v_i\in C$.  By assumption, 
we have the prox-regularity condition (2.2) holding for each $i$ for 
neighborhoods $\cX_i$ of $\bar x$, $\cV_i$ of $\bar v_i$, and parameter 
values $r_i$, $\alpha_i$.  Let $O$ denote an open convex neighborhood of 
the origin small enough that $\bar v_i+O\subset \cV_i$ for all $i$, and 
take $\cX=\cap_i\cX_{i}$, $\alpha=\min_i \alpha_{i}$, 
$r=\max_i r_{i}$.   Then $\cX$ is an open neighborhood of $\bar x$, 
and (2.2) holds for all $i$ with respect to $\cX$, $r$, $\alpha$, and the 
open neighborhood $\bar v_i+O$ of $\bar v_i$:
$$
  f(x')\geq f(x)+v_i\mdot[x'-x]-\frac{r}{2}|x'-x|^2 \; \forall x'\in\cX, \; 
    (x,v_i)\in(\cX^\alpha_f\times[\bar v_i+O]\,)\cap \gph\partial f.
$$
Consider now any $v\in \bar v +O$, writing it as $\bar v+w$ for some 
$w\in O$ and noting that  
$v=\sum_{i=0}^n\lambda_i v_i$ for $v_i=\bar v_i+w\in \bar v_i+O$.   From
the convex combination of the prox-regular inequalities for each $i$, we
see that the inequality likewise holds for $v$.    In particular, $v$ is
thus a proximal subgradient of $f$ at $\bar x$, so 
$v\in \partial f(\bar x)$.  This confirms that
$\conv C\subset\partial f(\bar x)$ and that (2.2) holds for $\cX$, $r$, 
$\alpha$ and the neighborhood $\cV=\bar v+O$ of $\bar v$, which gives 
the prox-regularity of $f$ at $\bar x$ for $\bar v$. 

A quantitative refinement of this argument, in relation to the level $r$,
will be important in (b) and (c). Refinement comes from the fact, coming
from (2.5), that the $r_{i}$-levels available for the prox-regularity 
at $\bar x$ for $\bar v_i$ include all the values in the interval
$(-\cnv f(\bar x\for\bar v_i), \infty)$.   It follows that, in taking
$r=\max_i r_{\bar v_i}$ we can arrange for $r$ to be any value
$>-\min_i \cnv f(\bar x\for\bar v_i)$.  Put another way, the argument
shows that, for any $\bar s\in\reals$,
$$
   \cnv f(\bar x\for\bar v)>\bar s,\,\forall \bar v\in C   \implies
   \cnv f(\bar x\for\bar v)>\bar s,\,\forall \bar v\in \conv C.
\eqno(2.14)$$ 
   
(b):   On the basis of (a) and the fact that the convex hull of a compact
set is compact, plus the observation in (2.14), we can henceforth assume 
for simplicity that $C$ itself is convex:  $C=\conv C$.  For each 
$\bar v\in C$ there are open neighborhoods $\cX_{\bar v}$ and 
$\cV_{\bar v}$ along with $\alpha_{\bar v} >f(\bar x)$ and $r_{\bar v}$ 
in the pattern of (2.2).  The union of the open sets $\cV_{\bar v}$ covers
 $C$, so the compactness of $C$ allows us to extract a finite covering, 
say for $\bar v_k$, $k=1,\ldots,m$.  Take $\cV=\cup_k\cV_{\bar v_k}$, 
$\cX=\cap_k\cX_{\bar v_k}$, $\alpha=\min_k \alpha_{\bar v_k}$ and 
$r=\max_k r_{\bar v_k}$.   Then $\cX$ is an open neighborhood of $\bar x$, 
$\cV$ is an open set containing $C$, and (2.2) holds.  

Here again we have a refinement in terms of the convexity modulus, because 
in arranging the covering we can insist, for any $\eps>0$, on having 
$r_{\bar v} < -\cnv f(\bar x\for\bar v)+\eps$.  That way, if $\bar s$ is a
value such that $\cnv f(\bar x\for\bar v)>\bar s$ for every $\bar v\in C$,
we can get $r<-\bar s$.
\ssk

(c): By choosing $\cX_{\bar v}$ and $\cV_{\bar v}$ in this construction to 
be convex, we get $\cX$ convex, but not necessarily $\cV$.  To fix that, 
replace $\cV$ by the open convex set $\lset v \mset\dist_C(v)<\eps\rset$ 
for a positive $\eps>0$ small enough that it lies within $\cV$.  Such 
$\eps$ must exist, since otherwise it would be possible to generate a 
$v$-sequence in the closed exterior of $\cV$ that gets closer and closer 
to $C$ and is bounded because $C$ is bounded.  It would have a cluster 
point in $C$ that still belongs to the exterior, in contradiction to 
$C\subset \cV$.  

The argument in (b) can now be sharpened by taking advantage of Theorem 2.1 
through its consequence in the three descriptions of the convexity modulus
in (2.5), (2.6) and (2.7).   By choosing smaller neighborhoods 
$\cX_{\bar v}$  and $\cV_{\bar v}$ and lower values for $\alpha_{\bar v}$ 
and $s_{\bar v}=-r_{\bar v}$, we can arrange to have not only the
prox-regularity condition (2.2), but also the variational convexity 
condition (2.3) for a function $\Hat f_{\bar v}$ and the associated 
maximal $r$-monotonicity in (2.4).  Proceeding with those elements to a
covering that corresponds to a finite collection of vectors $\bar v_k$, 
$k=1,\ldots,m$, we can pass again to $\cX$, $\cV$, $\alpha$ and $s$ as 
$-r$, know from the refinement at the end of the proof of (b) that if 
$\cnv f(\bar x,\bar v)>\bar s$ for all $\bar v\in C$, we can end up with 
$s>\bar s$.

Each of the $s$-convex functions $\Hat f_{\bar v_k}$ is defined on $\cX$ 
and $\leq f$ there, so by taking $\Hat f$ to be their pointwise max, we 
obtain an $s$-convex function that is $\leq f$ on $\cX$.   We must still 
confirm, though, that $\Hat f$ fits the rest of (2.3) as well.  For this, 
note that for any $k$ we have 
$$\eqalign{
  \Hat f_{\bar v_k}\leq f \text{on $\cX$,} \qd 
   (\cX\times[\cV\cap\cV_{\bar v_i}])\cap \gph\partial\Hat f_{\bar v_k} = 
   (\cX^\alpha_f\times[\cV\cap\cV_{\bar v_i}])\cap \gph\partial f &\cr
 \hbox{and, for all $(x,v)$ in this truncated graph, 
$\Hat f_{\bar v_k}(x)= f(x)$.}  &}
\eqno(2.15)$$ 
The truncation is furthermore maximal $s$-monotone relative to
$\cX\times[\cV\cap\cV_{\bar v_k}]$.  Note next that, because of common
$s$-convexity and $\Hat f_{\bar v_k} \leq \Hat f$, we have
$$
 \partial \Hat f_{\bar v_k}(x) \subset \partial \Hat f(x)
     \text{when} \Hat f_{\bar v_k}(x) =\Hat f(x).
$$ 
From this and having $\Hat f\leq f$ on $\cX$, we see in (2.3) that
$$
   (\cX\times[\cV\cap\cV_{\bar v_k}])\cap \gph\partial\Hat f_{\bar v_k}  
   \subset
   (\cX\times[\cV\cap\cV_{\bar v_k}])\cap \gph\partial\Hat f.
$$
Then equality must hold, however, because of the right side being
$s$-monotone and the left side being maximal $s$-monotone.  It follows
that (2.15) holds with $\Hat f$ in place of $\Hat f_{\bar v_k}$.  This
being true for every $k$, we conclude that (2.3) does hold for $\Hat f$.
\eop

A parametric property of partial subgradients will complete this section.
Several other special consequences of prox-regularity for subgradient 
calculus, such as a chain rule and a sum rule, were uncovered in 
\cite{CalcProxReg}, but the parameterization puts this new result in a 
different category.

\proclaim Theorem 2.7 {\rm (partial subgradients under prox-regularity)}. 
For closed proper $\phi$ on $\reals^n\times\reals^m$ and the mapping
$\,Y: (x,u,v)\mapsto \lset y \mset (v,y)\in\partial\phi(x,u)\rset$, suppose 
for given $(\bar x,\bar u,\bar v)$ that
  \paritem{(a)} for every $\bar y\in Y(\bar x,\bar u,\bar v)$,
prox-regularity of $\phi$ holds at $(\bar x,\bar u)$ for the subgradient 
$(\bar v,\bar y)$, and
  \paritem{(b)} the constraint qualification holds that
$$
     (0,y)\in\partial^\infty\phi(\bar x,\bar u) \text{only for} y=0.
\eqno(2.16)$$
Then there is a $\phi$-attentive neighborhood $(\cX\times\cU)_\phi^\alpha$ 
of $(\bar x,\bar u)$ along with open convex sets $\cV\ni \bar v$ and
$\cY\supset Y(\bar x,\bar u,\bar v)$ such that the mapping $Y$ is 
compact-convex-valued, outer semicontinuous and locally bounded from 
$(\cX\times\cU)_\phi^\alpha\times\cV$ into $\cY$, and
$$\eqalign{ 
  \hbox{$\forall u\in\cU,\,\exists\,x\in\cX$ with $\phi(x,u)<\alpha$ 
     and, for every such $x$,} \h{80} &\cr
  \partial_\low x\phi(x,u)\cap\cV =\lset v\in\cV\mset \exists\,y,\,
     (v,y)\in\partial\phi(x,u)\rset, \text{a convex set,} &}
\eqno(2.17)$$
and moreover such that uniform prox-regularity holds in the existence of 
$r>0$ for which
$$\eqalign{
  &\phi(x',u')\geq \phi(x,u)+(v,y)\mdot[(x',u')-(x,u)]
     -\frac{r}{2}|(x',u')-(x,u)|^2  \cr
  &\forall\,(x',u')\in\cX\times\cU, \text{when} (x,u,v,y)\in 
 [(\cX\times\cU)_\phi^\alpha\times(\cV\times\reals^n)]\cap\gph\partial\phi,}
\eqno(2.18)$$
with $r<\bar r$ if for every $\bar y$ in assumption (a) the prox-regularity
holds at a level $<\bar r$.

\state Proof.  The constraint qualification in (2.16) corresponds to the 
epigraphical mapping $u\mapsto \epi \phi(\cdot,u)$
having the Aubin property at $\bar u$ for the image point 
$(\bar x,\phi(\bar x,\bar u))$ \cite[10.16]{VA}.  That property includes
inner semicontinuity in the sense that for every neighborhood $\cW$ of
$(\bar x,\phi(\bar x,\bar u))$ the intersection of $\cW$ with 
$\epi\phi(\cdot,u)$ is nonempty for all $u$ near enough to $\bar u$.  This 
will later serve to justify the existence claim in (2.17).  But the Aubin
property persists by its definition to holding also for all $u$ near
enough to $\bar u$ with respect to $(x,\alpha)\in\epi\phi(\cdot,u)$ near
enough to $(\bar x,\phi(\bar x,\bar u))$.   Thus (2.16) persists
for every $(x,u)$ in some $\phi$-attentive neighborhood 
$(\cX_0\times\cU_0)_\phi^{\alpha_\low 0}$ of $(\bar x,\bar u)$.
We have then from \cite[10.11]{VA} that 
$$\eqalign{
  &\text{for} (x,u)\in (\cX_0\times\cU_0)_\phi^{\alpha_\low 0},\h{200} \cr
  & \partial_\low x\phi(x,u)\subset \lset v\mset 
   \exists\,y,\,(v,y)\in\partial\phi(x,u)\rset
   = \lset v\mset Y(x,u,v)\neq\emptyset\rset, }
\eqno(2.19)$$
with $Y(x,u,v)$ being compact.  Here, however, we need still better
understanding of how $Y(x,u,v)$ depends on $(x,u,v)$.  We can get that 
next by repeating the argument that's behind (2.19) and extracting from 
it the claimed outer semicontinuity and local boundedness of $Y$ at 
$(\tilde x,\tilde u,\tilde v)$  as a representative element of 
$((\cX_0\times\cU_0)_\phi^{\alpha_\low 0}\times\reals^n)\cap\dom Y$.  This 
means starting from
$$
    (v^\nu,y^\nu)\in\partial\phi(x^\nu,u^\nu) \text{with}
    (x^\nu,u^\nu,v^\nu)\to (\tilde x,\tilde u,\tilde v) \text{and}
      \phi(x^\nu,u^\nu)\to \phi(\tilde x,\tilde u),
\eqno(2.20)$$   
and demonstrating that the $y^\nu$ sequence must be bounded with all its 
cluster points $\tilde y$ satisfying 
$(\tilde v,\tilde y)\in\partial\phi(\tilde x,\tilde u)$.  That cluster point 
property is automatic from boundedness according to the definition of 
$\partial\phi$, so we can concentrate on establishing the boundedness.  

If boundedness failed, we could, by selection, reduce to having (2.20) with 
$0<|y^\nu|\upto\infty$.  Then, in taking $\lambda^\nu =1/|y^\nu|$ we could
further suppose that $\lambda^\nu y^\nu\to y$ with $|y|=1$.  In that case
$(\lambda^\nu v^\nu,\lambda^\nu y^\nu)\to (0,y)$ belonging to 
$\partial^\infty\phi(\tilde x,\tilde u)$ by the definition of that set
of horizon subgradients \cite[8.3]{VA}.  But that's precluded by our 
assumption (b).  This argument confirms not only (2.19) for a neighborhood 
$(\cX_0\times\cU_0)_\phi^{\alpha_\low 0}$, but also the outer semicontinuity 
and local boundedness of $Y$ on 
$(\cX_0\times\cU_0)_\phi^{\alpha_\low 0}\times\reals^n$.

The next step is applying Proposition 2.6 to $\phi$ and the compact set 
$C=\lset (\bar v,\bar y) \mset \bar y\in Y(\bar x,\bar u,\bar v)\rset$
within $\partial\phi(\bar x,\bar u)$.   From Proposition 2.6(a) we get
$\conv C\subset\partial\phi(\bar x,\bar u)$, hence 
$\conv Y(\bar x,\bar u,\bar v) \subset Y(\bar x,\bar u,\bar v)$, i.e., 
the convexity of $Y(\bar x,\bar u,\bar v)$.   From Proposition 2.6(b) and 
(c) we get open convex sets $\cW\ni (\bar x,\bar u)$ and $\cZ\supset C$ 
along with $\alpha>\phi(\bar x,\bar u)$ and $r$ for which
$$\eqalign{
  &\phi(x',u')\geq \phi(x,u)+(v,y)\mdot[(x',u')-(x,u)]
     -\frac{r}{2}|(x',u')-(x,u)|^2  \cr
  &\hbox{for all}\;\, (x',u')\in\cW \text{when} (x,u,v,y)\in 
   (\cW_\phi^\alpha\times\cZ)\cap\gph\partial\phi,  }
\eqno(2.21)$$
with $r<\bar r$ if the prox-regularity in assumption (a) is $<\bar r$
for every $\bar y$.

We can think of $C$ as $Y_\pls(\bar x,\bar u,\bar v)$ for the mapping
$Y_\pls(x,u,v) = \{v\}\times Y(x,u,v)$, which inherits from $Y$ its outer 
semicontinuity and local boundedness relative to 
$(\cX_0\times\cU_0)_\phi^{\alpha_\low 0}\times\reals^n$ for the neighborhood 
$(\cX_0\times\cU_0)_\phi^{\alpha_\low 0}$ in (2.19).  For $(x,u,v)$ near 
enough to $(\bar x,\bar u,\bar v)$ there, we'll have 
$Y_\pls(x,u,v)\subset\cZ$.  That will yield the convexity of 
$Y_\pls(x,u,v) = \cZ \cap \partial\phi(x,u) \cap (\{v\}\times \reals^m)$
by Proposition 2.5(b) and with it the convexity of $Y(x,u,v)$.
Without loss of generality we can have the open convex neighborhoods
in (2.21) be small enough for that to hold and also be of the form
$$
  \cW_\phi^\alpha=(\cX\times\cU)_\phi^\alpha
  \subset (\cX_0\times\cU_0)_\phi^{\alpha_\low 0} \qquad \cZ=\cV\times\cY,
\eqno(2.22)$$
with $\cX$, $\cU$, $\cV$, $\cY$, all convex.  That gives us (2.18).  In the 
case of (2.18) with $u'=u$, we have 
$$\eqalign{
   &\phi(x',u)\geq \phi(x,u)+v\mdot[x'-x]-\frac{r}{2}|x'-x|^2, 
    \,\forall x'\in\cX, \cr
   &\text{when}    
    (x,u)\in (\cX\times\cU)_\phi^\alpha,\;\, v\in\cV,\;\, y\in Y(x,u,v), }
$$
which implies $v\in\partial_x\phi(x,u)\cap\cV$.   This verifies that (2.19) 
for $(\cX\times\cU)_\phi^\alpha$ must hold with equality, hence (2.17).  
The convexity of $\partial_x\phi(x,u)\cap\cV$ follows then from 
Proposition 2.5(b) and the convexity of $\cV$, inasmuch as $\phi(\cdot,u)$ 
inherits prox-regularity from $\phi$ under the subgradient formula in 
(2.17).  \eop

\proclaim Corollary 2.8 {\rm (uniform prox-regularity of partial functions)}.
For $\phi$ as in Theorem 2.7, there exist open convex neighborhoods $\cX$ 
of $\bar x$, $\cV$ of $\bar v$ and $\cU$ of $\bar u$ along with 
$\alpha> \phi(\bar x,\bar u)$ and $r$ such that, for each $u\in\cU$, the 
partial function $\phi^u(x)=\phi(x,u)$ satisfies the prox-regularity 
condition
$$\eqalign{
    \phi^u(x') \geq \phi^u(x)+v\mdot[x'-x]-\frac{r}{2}|x'-x|^2, \;
      \forall x'\in\cX,  &\cr
 \hbox{when}\;\, (x,v)\in(\cX\times\cV)\cap \gph\partial \phi^u 
      \text{with} \phi^u(x)<\alpha, &}
\eqno(2.23)$$
moreover with $r<\bar r$ if for every $\bar y\in Y(\bar x,\bar u,\bar v)$
the assumed prox-regularity holds at a level $<\bar r$.

\state Proof.   This comes out of (2.17) and (2.18), inasmuch as 
$\partial \phi^u(x)=\partial_x\phi(x,u)$.   \eop

\section{
                Main proofs and counterexamples
 } 

This section takes up the task of proving the theorems stated in Section 1  
with the help of the background and additional results in Section 2, while 
also offering further insights through examples.  A closer look at
variational sufficiency comes first.

Recall that the strong variational sufficient condition for local 
optimality in problem $(\bar P)$ is fulfilled at $\bar x$ in partnership 
with a multiplier vector $\bar y\in \Bar Y =\lset \bar y \mset
(0,\bar y)\in\partial\phi(\bar x,0)\rset$ if there is an elicitation level 
$e\geq 0$ such that the function 
$$
   \phi_e(x,u)=\phi(x,u)+\frac{e}{2}|u|^2,  \text{with}
     \partial \phi_e(x,u)=\partial\phi(x,u)+(0,eu),
\eqno(3.1)$$ 
is variationally $s$-convex at $(\bar x,0)$ for $(0,\bar y)$ at a level 
$s>0$.  In the terminology (1.7), strong variational suffiency at $\bar x$
without mention of a particular $\bar y$ refers to having this condition
fulfilled for every $\bar y\in\Bar Y$, but conceivably this could be with 
different levels of $e$ and $s$ that depend on $\bar y$.  In fact, as
demonstrated next, the same $e$ and $s$ can be made to work, even with the
same neighborhoods.

\proclaim Proposition 3.1 {\rm (uniformity in strong variational sufficiency)}.
Strong variational sufficiency for local optimality at $\bar x$ under the
basic constraint qualification (1.5) corresponds to having open convex sets 
$\cX\times\cU\ni (\bar x,0)$ and $\cV\times\cY \supset \{0\}\times\Bar Y$, 
parameter levels $e\geq 0$, $s>0$, $\alpha>\phi(\bar x,0)$, and a closed 
proper $s$-convex function $\psi\leq\phi_e$ on $\cX\times\cU$ such that
$$\eqalign{
 (x,u,v,y)\in \gph\partial\psi\cap[(\cX\times\cU)\times(\cV\times\cY)] &\cr
 \iff (x,u,v,y)\in \gph\partial\phi_e \cap 
    [(\cX\times\cU)^\alpha_{\phi_e}\times(\cV\times\cY)], &\cr
   \hbox{and $\psi(x,u) = \phi_e(x,u)$ 
   for such vectors $(x,u,v,y)$,} & } 
\eqno(3.2)$$
with $\phi$ then having the properties in Theorem 2.7 for $\bar u=0$, 
$\bar v=0$, and $r=-s$.  This furthermore entails uniform variational 
$s$-convexity of the partial functions $\phi^u(x)=\phi(x,u)$ in the sense 
of having for all $u\in\cU$, with respect to 
$\psi^{u,e}(x):=\psi(x,u)-\frac{e}{2}|u|^2$, for which
$\psi^{u,e}\leq\phi^u$ on $\cX$,
$$\eqalign{
 (x,v)\in \gph\partial\psi^{u,e}\cap(\cX\times\cV) &\cr
 \iff (x,v)\in \gph\partial\phi^u \cap 
    (\cX^\alpha_{\phi^u}\times\cV), &\cr
   \hbox{and $\psi^{u,e}(x) = \phi^u(x)$ for such vectors $(x,v)$.} & } 
\eqno(3.3)$$
        
\state Proof.  For each $\bar y\in\Bar Y$ we are given the existence of 
open neighborhoods $\cX_{\bar y}\times\cU_{\bar y}$ of $(\bar x,0)$ and
$\cV_{\bar y}\times\cY_{\bar y}$ of $(0,\bar y)$, parameter values
$e_{\bar y}\geq 0$, $s_{\bar y}>0$, $\alpha_{\bar y}>\phi(\bar x,0)$, and
a closed proper $s_{\bar y}$-convex function $\psi_{\bar y}$  on
$\cX_{\bar y}\times\cU_{\bar y}$ such that 
$$\eqalign{
   (x,u,v,y)\in \gph\partial\psi_{\bar y}\cap
  [(\cX_{\bar y}\times\cU_{\bar y})\times(\cV_{\bar y}\times\cY_{\bar y})] 
               &\cr
 \iff (x,u,v,y)\in \gph\partial\phi_{e_{\bar y}} \cap 
    [(\cX_{\bar y}\times\cU_{\bar y})^{\alpha_{\bar y}}_{\phi_{e_{\bar y}}}
      \times(\cV_{\bar y}\times\cY_{\bar y})], &\cr
   \hbox{and $\psi_{\bar y}(x,u) = \phi_{e_{\bar y}}(x,u)$ 
                for such vectors $(x,u,v,y)$.} &} 
\eqno(3.4)$$
The key observation is that these same elements serve equally then to
describe the strong variational sufficiency at $\bar x$ not only for 
$\bar y$, but also for every $\bar y'\in\Bar Y\cap\cY_{\bar y}$.  The 
compactness of $\Bar Y$, assured by the constraint qualification,
allows us to select finitely many $\bar y_k\in\Bar Y$ such that the open
neighborhoods $\cY_{\bar y_i}$ cover $\Bar Y$.  Then for 
$e=\max_k e_{\bar y_k}\geq 0$ and $s=\min_k s_{\bar y_k}>0$ we have, for 
every $\bar y\in\Bar Y$, that $\phi_e$ is variationally $s$-convex at
$(\bar x,0)$ with respect to $(0,\bar y)$.  In particular, for any 
$\bar s\in (0,s)$, we have $\cnv\phi_e(\bar x,0\for 0,\bar y) >\bar s$ for
every $(0,\bar y)\in C=\{0\}\times\Bar Y =\partial\phi(\bar x,0)$ and are
able that way to apply Proposition 2.6(c) to reach the uniform condition
in (3.2), with a slightly lower $s>0$ if necessary.

Since prox-regularity at level $r=-s$ is implied by such variational 
$s$-convexity, Theorem 2.7 is applicable, as claimed, for the case of
$\bar u=0$, $\bar v=0$, and $\alpha>\phi_e(\bar x,0)\geq\phi(\bar x,0)$,
that carries over to $\phi_e$ in yielding
$$\eqalign{ 
  \hbox{$\forall u\in\cU,\,\exists\,x\in\cX$ with $\phi_e(x,u)<\alpha$ 
     and, for every such $x$,} \h{80} &\cr
  \partial_\low x\phi_e(x,u)\cap\cV =\lset v\in\cV\mset \exists\,y,\,
     (v,y)\in\partial\phi_e(x,u)\rset. &}
\eqno(3.5)$$
We likewise have for $(x,u)\in\cX\times\cU$ that $\partial_x \psi(x,u)=
\lset v\in\cV\mset \exists y,\,(v,y)\in\partial\psi(x,u)\rset$ if the
constraint qualification holds that $(0,y)\in\partial^\infty\psi(x,u)$
only for $y=0$.  That does hold because the convexity of $\psi$ reduces 
it to $u$ belonging to 
$\,\iint \lset u \mset \exists\,x,\,\psi(x,u)<\infty\rset$, which includes 
$\cU$ through (3.5) since here $\psi(x,u)\leq\phi_e(x,u)$.  These
observations about partial subgradients, in combination with noting that
$\psi(x,u)=\phi_e(x,u)$ corresponds to $\psi^{u,e}(x)=\phi^u(x)$,
translate (3.2) into (3.3) for fixed $u$.  \eop

\state Proof of Theorem 1{.}1.  In all three conditions (A), (B) and 
(C), we have the basic constraint qualification (1.5) combined with 
$\phi$ being prox-regular at $(\bar x,0)$ for every $(0,\bar y)$ with 
$\bar y\in\Bar Y$, so we are always in the setting of Theorem 2.7 and  
have at our disposal all the subgradient and prox-regularity properties 
of $\phi$ there for $\bar v=0$ and $\bar u=0$.  

(A)$\Rightarrow$(B):  Along with the prox-regularity of $\phi$ at
$(\bar x,0)$ for every $(0,\bar y) \in\partial\phi(\bar x,0)$, i.e., 
$\bar y\in\Bar Y$, we have the Lipschitz continuity of $M_\delta(v,u)$ with 
respect to $(v,u)\in\cV\times\cU$ for a Lipschitz constant $\kappa$.  Fix 
any $\bar y\in\Bar Y$.  To reach (B), it must be shown that, for $e>0$ 
high enough, the function $\phi_e(x,u)=\phi(x,u)+\frac{e}{2}|u|^2$ is 
variationally strongly convex at $(\bar x,0)$ for the subgradient 
$(0,\bar y)$.  By Theorem 2.2 as applied to $\phi_e$ in place of $f$, 
that's equivalent to the tilt stability of the function
$$
   \Bar\phi_e(x,u):= \phi(x,u)-\bar y\mdot u+\frac{e}{2}|u|^2 
\eqno(3.6)$$ 
at $(\bar x,0)$.  According to Theorem 2.4, we can therefore confirm it by 
demonstrating that $\Bar\phi_e$ 
has a local minimum at $(\bar x,0)$ and satisfies
$$
   (0,0)\in\partial_*^2\Bar\phi_e(\bar x,0;0,0)(\xi,\omega)
     \text{only for} (\xi,\omega)=(0,0).
\eqno(3.7)$$
Here $\partial\Bar\phi_e(x,u)=\partial\phi(x,u)+(0,eu-\bar y)$ and
$\partial_*^2\Bar\phi_e(\bar x,0;0,0)(\xi,\omega) =
   \partial_*^2\phi(\bar x,0;0,0)(\xi,\omega) +(0,e\omega)$, so that showing
(3.7) comes down to showing
$$ 
   (0,-e\omega)\in\partial_*^2\phi_e(\bar x,0;0,0)(\xi,\omega)
     \text{only for} (\xi,\omega)=(0,0).
\eqno(3.8)$$

To establish the local minimum for $e$ sufficiently high, we invoke the
prox-regularity in (2.18) of Theorem 2.7 for $\bar v=0$ and $\bar u=0$ to
see first that
$$\eqalign{
    \phi(M_\delta(0,u),u)\geq \phi(\bar x,0)+ 
       (0,\bar y)\mdot[(M_\delta(0,u),u)-(\bar x,0)]
       -\frac{r}{2}|(M_\delta(0,u),u)-(\bar x,0)|^2,  &\cr
   \text{where} \phi(x,u) \geq \phi(M_\delta(0,u),u) 
        \text{and} \bar x=M_\delta(0,0),
     \text{hence} |M_\delta(0,u)-\bar x|\leq\kappa|u|.  &}
$$
We get from this that
$$
   \Bar\phi_e(x,u)-\Bar\phi_e(\bar x,0)\geq
     \frac{e}{2}|u|^2-\frac{r}{2}(\kappa^2|u|^2+|u|^2)
      =\frac{e-r(\kappa^2+1)}{2}|u|^2.
$$
To ensure the local minimum, it suffices therefore to have
$$
      e > r(\kappa^2+1).
\eqno(3.9)$$

Turning now to verifying the strict second-order subdifferential condition 
in (3.8), recall from the definition of that mapping in (2.12), in its
articulation for $\phi$, that
$$\eqalign{
 (\xi,\omega,0,-e\omega)\in\gph \partial^2_*\phi(\bar x,0\for 0,\bar y)
          \iff &\cr
\h{40}     \exists\, (x^\nu,u^\nu,v^\nu,y^\nu)\to(\bar x,0,0,\bar y) 
         \txt{in}\gph\partial\phi &\cr
\h{40}  \hbox{and}\;\, 
    (\xi^\nu,\omega^\nu,\mu^\nu,\eta^\nu)\to(\xi,\omega,0,-e\omega),\; 
        \tau^\nu\dnto 0, \txt{with}  &\cr
\h{40}    (x^\nu+\tau^\nu\xi^\nu,u^\nu+\tau^\nu\omega^\nu,
  v^\nu+\tau^\nu\mu^\nu,y^\nu+\tau^\nu\eta^\nu)\in\gph\partial\phi, 
     \txt{while} &\cr
\h{40} \phi(x^\nu,u^\nu)\to \phi(\bar x,0) \text{and}
  \phi(x^\nu+\tau^\nu\xi^\nu,u^\nu+\tau^\nu\omega^\nu)\to\phi(\bar x,0). &}
\eqno(3.10)$$
In these circumstances we have $x^\nu=M_\delta(v^\nu,u^\nu)$ and 
$x^\nu+\tau^\nu\xi^\nu=
M_\delta(v^\nu+\tau^\nu\mu^\nu,u^\nu+\tau^\nu\omega^\nu)$, hence
$$
   |\xi^\nu|=\frac{1}{\tau^\nu}
    \Big|M_\delta(v^\nu+\tau^\nu,u^\nu+\tau^\nu\omega^\nu)
  -M_\delta(v^\nu,u^\nu)\Big| \leq \kappa|(\mu^\nu,\omega^\nu)|,
$$
and in the limit
$$
      |\xi|\leq \kappa|\omega|.
\eqno(3.11)$$      
On the other hand, the prox-regularity coming from Theorem 2.7 with $r>0$
in (2.18) gives us from (3.10), through the associated monotonicitiy of 
$\partial\phi$ at level $s=-r$ in Theorem 2.1, that
$$\eqalign{
 \Big((v^\nu+\tau^\nu\mu^\nu,y^\nu+\tau^\nu\eta^\nu)-(v^\nu,y^\nu)\Big)\mdot
 \Big((x^\nu+\tau^\nu\xi^\nu,u^\nu+\tau^\nu\omega^\nu)-(x^\nu,u^\nu)\Big)
   &\cr
 \h{30} \geq 
 -r|(x^\nu+\tau^\nu\xi^\nu,u^\nu+\tau^\nu\omega^\nu)-(x^\nu,u^\nu)|^2. &}
$$
This reduces to
$
 (\mu^\nu,\eta^\nu)\mdot(\xi^\nu,\omega^\nu)\geq -r|(\xi^\nu,\omega^\nu)|^2,
$
which in passing to the limit in (3.10) turns into
$(0,-e\omega)\mdot(\xi,\omega)\geq -r|(\xi,\omega)|^2$, or in other words
$$
       r|\xi|^2 \geq (e-r)|\omega|^2.
\eqno(3.12)$$
In comparing this with (3.11), which entails $r|\xi|^2\leq 
r\kappa^2|\omega|^2$, we obtain $0\geq [e-r(\kappa^2+1)]|\omega|^2$.  The
conclusion is that, as long as $e$ is large enough to satisfy (3.9), we get
from (3.12) that $\omega=0$ and then from (3.11) that also $\xi=0$, so that
(3.8) is correct.

(B)$\Rightarrow$(C):   Proposition 3.1 furnishes in (3.3) variational 
strong convexity in $x$ component that through Theorems 2.2 and 2.3(c) 
produces, for every $u\in\cU$, tilt stability in $v\in\cV$ with Lipschitz 
constant $s^{-1}$.   That aspect of (C) in turn underpins, through 
Theorem 2.3(a)(b), to $M_\delta$ being the localization 
$M^{\delta,\alpha}$ of $M$ in (1.8) for sufficiently low $\alpha$, and 
$m_\delta(v,u)$ having the properties with respect to $v$ 
that are indicated in the theorem.  

For the continuity of $M_\delta(v,u)$ in $u$ and the additional properties 
of $m_\delta(v,u)$ with respect to $u$ claimed under (B), further argument 
is needed.  Consider the replacement of $\phi$ by $\phi_e$ in the formulas 
(1.1) for $M_\delta(v,u)$ and $m_\delta(v,u)$.  This would have no effect 
on $M_\delta(v,u)$ but would add $\frac{e}{2}|u|^2$ to $m_\delta(v,u)$; let 
that modified function be denoted by $m^e_\delta(v,u)$.   In those terms 
we have, for $u\in\cU$,
$$
  m^e_\delta(v,u)=\min_{|x-\bar x|<\delta}
      \lset\phi_e(x,u)-v\mdot[x-\bar x]\,\rset,
          \text{attained for} x=M_\delta(v,u)=M^{\delta,\alpha}(v,u).
$$
Then, from the definition of $M^{\delta,\alpha}(v,u)$ in (1.8) and via 
(1.3), we have some $y$ such that $(v,y)\in\partial\phi(x,u)$, hence 
$(v,y_e)\in\partial\phi_e(x,u)$ for $y_e=y+eu$, all under the umbrella of 
(3.2), so that actually $\phi_e(x,u)=\psi(x,u)$.   This reveals that
$$
  m^e_\delta(v,u)=\min_{|x-\bar x|<\delta}
         \lset \psi(x,u)-v\mdot[x -\bar x]\, \rset,
          \text{attained for} x=M_\delta(v,u).
\eqno(3.13)$$
Since $\psi$ is in particular convex, this inf-projection formula
implies that $m^e_\delta(v,u)$ is convex in $u\in \cU$ and moreover by
duality \cite[11.39]{VA} that 
$
   \partial_u m^e_\delta(v,u)= Y_e(x,u,v)=Y(x,u,v)+eu 
            \text{for} x=M_\delta(v,u).
$
Translating this back to $m_\delta(v,u)$, we get the claimed formula (1.9).  

The continuity of $M_\delta(v,u)$ in $u$ (which becomes continuity in
$(v,u)$ when combined with the Lipschitz continuity in $v$) also comes out 
of (3.13), as follows.  Consider the strongly convex function $f_v(x,u)$ 
defined by
$$
 f_v(x,u)= \psi(x,u)-v\mdot[x -\bar x] \text{when}
    |x-\bar x|\leq\delta, \text{otherwise} f_v(x,u)=\infty.
$$
Observe from (3.13) and the strong convexity that the singleton 
$x=M_\delta(v,u)$ is $\argmin f_v(\cdot,u)$.  According to \cite[7.44]{VA} 
in describing parametric minimization under convexity, the mapping 
$u\mapsto\argmin f_v(x,u)$ is continuous on the interior of the convex set 
$U$ consisting of the vectors $u$ such that $f_v(\cdot,u)$ is proper.  
That set $U$ obviously includes the vectors $u$ where $\argmin f_v(\cdot,u)$ 
is nonempty, so its interior includes the open neighborhood $\cU$ of $0$ 
where we have singleton $M_\delta(v,u)$.  \eop

The example presented next demonstrates that the strong variational
sufficiency in (B) is inadequate for having the full stability in (A),
even when $\phi$ itself is already a strongly convex function, in which
case no elicitation is even needed.

\proclaim Example 3.2 {\rm (full stability absent despite strong convexity)}.
For $(x,u)\in\reals\times\reals$ and $(\bar x,\bar u)=(0,0)$, the function
$\,\phi(x,u) = \frac{3}{4}|(x,u)|^\frac{4}{3} +|(x,u)|-x\,$
is strongly convex on any bounded convex neighborhood of $(0,0)$, where it
achieves its global minimum value of 0.  In particular it is prox-regular 
globally with the basic constraint qualification (1.5) satisfied.  However,
in the notation (1.1) for any $\delta>0$, it has $M_\delta(v,u)$ being
locally Lipschitz continuous in $v$ but not locally Lipschitz continuous 
in $(v,u)$.  Full stability is thus lacking.

\state Detail.  For a convex function $\phi$, $\partial^\infty\phi(x,u)$ 
is the normal cone to the effective domain of $\phi$ at $(x,u)$.  Here that
domain is $\reals\times\reals$, so the basic constraint qualification (1.5) 
is surely satisfied.  Local strong convexity of $\phi$ is the property dual 
to its conjugate $\phi^*$ being differentiable with $\nabla\phi$ locally
Lipschitz continuous, cf.\ \cite[Proposition 3.5]{ProxRegTests}, and here 
the conjugate function is $\phi^*(v,y)=\frac{1}{4}\dist^4(v,y\mid D)$, 
where $D$ = unit disk around $(1,0)$, which does fit that prescription.  
The strong convexity in $x$ for fixed $u$ yields global tilt stability and 
with it the single-valuedness of $M_\delta(v,u)$ with $M_\delta(0,0)=0$, 
where $M_\delta$ can be replaced by $M$ through globality.  But it 
will be seen that $M(0,u)$ fails to be Lipschitz continuous around 
$u=0$.  At $(v,u)$ apart from $(0,0)$, $\phi$ is continuously differentiable
with $x=M(v,u)$ if and only if $\frac{\partial\phi}{\partial x}(x,u)=v$.
In terms of the auxiliary variable $z=\sqrt{x^2+u^2} =|(x,u)|$, we have
$$
    \phi(x,u) = k(z)-x \text{for} k(z)=\frac{3}{4}z^\frac{4}{3}+z,
$$
hence
$
  \frac{\partial\phi}{\partial x}(x,u)=k'(z)\frac{\partial z}{\partial x}-1
     = \big(z^\frac{1}{3} +1\big)(x/z)-1,
$
so that $x=M(0,u) \Leftrightarrow x=z/\big(1+z^\frac{1}{3}\big)$.
From this and the definition of $z$ we get
$u^2=z^2-x^2 = z^2 - \big(z^2/\big(1+z^\frac{1}{3})^2\big) 
  = \big[z^2/\big(1+z^\frac{1}{3}\big)^2\big]
        \big[\big(1+z^\frac{1}{3}\big)^2-1\big]
     = x^2 \big[\big(1+z^\frac{1}{3}\big)^2-1\big]$.
Therefore $u^2/x^2 =\big(1+z^\frac{1}{3}\big)^2-1 \to 0$ as 
$(x,u)\to (0,0)$, which is equivalent to $M(0,u)/u \to \infty$ as 
$u\to 0$, where $M(0,0)=0$.  Thus, $M(0,u)$ lacks Lipschitz continuity 
around $u=0$, as claimed.  \eop

It might be wondered whether it's really necessary in (B) to require
strong variational sufficiency at $\bar x$ to hold for all multiplier
vectors $\bar y$.  Maybe if it holds for one it automatically holds for
all?   Here is a counterexample to that in the context of classical
nonlinear programming.

\proclaim Example 3.3 {\rm (strong variational sufficiency for some
multipliers but not all)}.
In the case of $\phi$ in $(P)$ that corresponds to the problem of 
minimizing $f_0(x_1,x_2,x_3,x_4)$ subject to 
$f_i(x_1,x_2,x_3,x_4)\leq 0$ for $i=1,2,3,4$, where
$$\eqalign{
   f_0(x_1,x_2,x_3,x_4) = x_3 + \half x_4^2   &\cr
   f_1(x_1,x_2,x_3,x_4) = x_1 -x_3   &\cr
   f_2(x_1,x_2,x_3,x_4) = -x_1 -x_3   &\cr
   f_3(x_1,x_2,x_3,x_4) = x_2 -x_3 - \half x_4^2   &\cr
   f_4(x_1,x_2,x_3,x_4) = -x_2 -x_3 - \half x_4^2,   &}
$$
with canonical perturbation parameters $u_i$, namely
$$  
   \phi(x,u)=f_0(x)+\delta_\low{\reals^4_\mns}
     (f_1(x)+u_1,f_2(x)+u_2,f_3(x)+u_3,f_4(x)+u_4), 
$$
the strong variational sufficient condition is satisfied at 
$\bar x=(0,0,0,0)$ for just some of the associated multiplier vectors
$\bar y$.

\state Detail.  The constraints require $x_3\geq |x_1|$ and 
$x_3\geq |x_2|-\half x_4^2$, so it's clear that $\bar x=(0,0,0,0)$ is 
locally optimal with all four constraints being active.  The subspace 
orthogonal to the four active constraint gradients is
$
      S = (0,0,0,\reals),
$
while the Lagrangian function is
$$\eqalign{
   L(x,y)= x_3+\half x_4^2 +y_1(x_1-x_3)-y_2(x_1+x_3)
                 +y_3(x_2-x_3- \half x_4^2)-y_4(x_2+x_3+ \half x_4^2)   &\cr
   \h{36}  = x_1(y_1-y_2) +x_2(y_3-y_4) +x_3(1-y_1-y_2-y_3-y_4)
    +\half(1-y_3-y_4)x_4^2.   &}
$$
The KKT conditions require $y_i\geq 0$ along with $y_1=y_2$, $y_3=y_4$, 
and $y_1+y_2+y_3+y_4=1$.  There is just one degree of freedom, so the set 
of multiplier vectors $\bar y$ associated with $\bar x$ can be
described as
$$
    \Bar Y= \lset \half(1-\theta,1-\theta,\theta,\theta) 
                      \mset 0\leq\theta\leq 1\rset.
$$
The classical strong second-order sufficient condition for local 
optimality, which is known from \cite[Example 1]{HiddenConvexity} to be 
equivalent to strong variational sufficiency, is satisfied at $\bar x$ 
for $\bar y\in Y$ if and only if the hessian 
$H=\nabla^2_{xx}L(\bar x,\bar y(\theta))$ is positive-definite relative 
to $S=(0,0,0,\reals)$.  The $4\times 4$ matrix $H$ has entries $h_{ij}$ that 
are 0 except possibly for $h_{44}$, with the strong second-order sufficient
condition be satisfied if and only if $h_{44}>0$.   But 
$h_{44}= 1 -2\theta$.  Thus we have the condition satisfied for the 
multiplier vectors $\bar y$ corresponding to $\theta\in [0,\half)$, but not 
for the ones corresponding to $\theta\in[\half,1]$.  \eop
\msk

\state Proof of Theorem 1{.}2.    From (C) we have the uniform Lipschitz continuity of $M_\delta(v,u)$ in $v\in\cV$ for each
$u\in\cU$.  To reach (A), we also need, for a possibly smaller $\delta$
and neighborhood $\cV\times\cU$ of $(0,0)$, the uniform Lipschitz 
continuity of $M_\delta(v,u)$ in $u\in\cU$ for each $v\in\cV$.  Here we 
can replace $M_\delta(v,u)$ by $M^{\delta,\alpha}(v,u)$ on the basis of 
the relationship (1.8) in Theorem 1.1.  The focus then is on having the 
mappings $M_v^{\delta,\alpha}$ in (1.10) be Lipschitz continuous around 
$u=0$ uniformly with respect to $v$ in some neighborhood of $v=0$.

Each mapping $M_v^{\delta,\alpha}$ is already single-valued and continuous 
around $u=0$ for $v$ near to $0$, so Lipschitz continuity is identical here 
to the Aubin property, although that ordinarily would be weaker.  A criterion 
for $M_v^\delta$ to have the Aubin property at $u=0$ with respect to 
$x=M_v^{\delta,\alpha}(0)$ is available from \cite[Theorem 4B.2]{DontRock}; 
it supposes $\gph M_v^{\delta,\alpha}$ to be locally closed around $(u,x)$ 
for $x= M_v^{\delta,\alpha}(u)$, but we know that does hold under the 
prox-regularity, because the optimal value 
$m_\delta (v,u)=\phi(M_v^{\delta,\alpha}(u),u)$ approaches 
$m_\delta (v,0)=\phi(M_v^{\delta,\alpha}(0),0)$ as $u\to 0$.  The criterion
in question is the uniform boundedness of the inner norms 
$|DM_v^{\delta,\alpha}(u)|^\mns$ on some neighborhood of $u=0$; such a 
bound then provides a Lipschitz constant.  Here we furthermore want 
uniformity with respect to $v$ in a neighborhood of $v=0$.  The finiteness 
of the upper limit in (1.12) guarantees all of that.   \eop



\msk 
\noindent
{\bf Data availability:   } no data was generated or analyzed.

\end{document}